\documentclass[12pt]{article}
\usepackage{amsfonts}\usepackage{hyperref}
\usepackage{pdfpages}
\usepackage{authblk}
\usepackage{graphicx}
\usepackage{amsthm}
\usepackage{amsmath}
\usepackage{amssymb}
\usepackage{fancyhdr}
\usepackage{parskip}
\usepackage{tikz}
\usepackage{hyperref}
\usepackage{setspace}
\usepackage[T1]{fontenc}
\usepackage{babel}
\usepackage{authblk}

\hypersetup{colorlinks=true,
urlcolor=blue,}

\usetikzlibrary{graphs}
\usetikzlibrary{graphs.standard}

\theoremstyle{plain}
\newtheorem{theorem}{Theorem}[section]

\newtheorem{cor}[theorem]{Corollary}
\newtheorem{lemma}[theorem]{Lemma}

\newtheorem{conjecture}[theorem]{Conjecture}

\newtheorem*{theorem*}{Theorem}
\newtheorem*{definition*}{Definition}
\newtheorem*{prop*}{Proposition}
\newtheorem*{lemma*}{Lemma}
\newtheorem*{corollary*}{Corollary}
\newtheorem*{conjecture*}{Conjecture}

\begin{document}
\title{Acyclic orientations and poly-Bernoulli numbers}

\author[1]{Peter~J.~Cameron}
\author[2]{ C.~A.~Glass}
\author[3]{Kamilla Rekv\'enyi}
\author[4]{R.~U.~Schumacher}
\affil[1]{University of St Andrews}
\affil[2]{City,University of London}
\affil[3]{Imperial College London}
\affil[4]{City,University of London}

\date{}
\maketitle

\begin{abstract}

The main contribution of this paper is a formula for the number of acyclic orientations of a complete bipartite, $K_{n_1,n_2},$ revealing that it is equal to the poly-Bernoulli number $B_{n_1}^{(-n_2)}$ introduced in 1997 by  Kaneko. We also give a simple bijective identification of acyclic orientations and lonesum matrices, confirming a 2008 result of Brewbaker, and show that the poly-Bernoulli numbers behave concavely, as well being symmetric.

A second goal is to explore the behaviour of more general complete $r$-partite graphs in the space of the number of acyclic orientations. 
We prove that the number of acyclic orientations of complete bipartite graphs on $n$ vertices is a unimodular concave function maximised by the Tur\'an graph $T(2,n).$ 
For tripartite graphs,  we derive an explicit formula for its number of its acyclic orientations, for the specific case with a single vertex partition. For more complex complete $r$-partite graphs an algorithmic approach is suggested, taking advantage of the relationship with graph colouring through the result of Stanley, namely the value of the graph colouring polynomial at the value -1. 

An underlying theme is the exploration of the space of the number of acyclic orientations of graphs. We present  various theoretical and computational results, with additional, conjectures for both complete $r$-partite and general graphs to guide future research.

\end{abstract}

\section{Introduction}

An \emph{acyclic orientation} of a graph $G$ is an assignment of direction to each edge in such a way that we obtain no directed cycles, thus obtaining an acyclic directed graph. 
Let $a(G)$ be the number of acyclic orientations of a graph
$G$. 
There is always at least one acyclic orientation of $G$ obtained by ordering the vertices of $G$ and orienting edges from smaller to larger vertex index.

There are connections between acyclic orientations and proper vertex colourings of $G$. For example, given a vertex colouring, we can obtain an acyclic orientation by ordering the colour classes and directing edges from smaller to greater.
Conversely, given an acyclic orientation, we can obtain a colouring by giving the first colour to the vertices with in-degree zero, then recursively colouring the remainder of the graph with the other colours.
In addition, Stanley~\cite{stanley} showed that, if $P_G$ denotes the chromatic polynomial of $G$ (so that $P_G(q)$ is the number of colourings with $q$ colours for positive integers $q$), then $a(G)=P_G(-1)$.

In the fourth author's PhD thesis, the following question was considered:
\begin{quote}
Given $m$ and $n$, what are the minimum and maximum values of $a(G)$
over all graphs with $n$ vertices and $m$ edges?
\end{quote}
The minimum value is known exactly, and the graphs realising it are completely characterised \cite{minimal}. 
However, the maximum is a different matter. 
The Hanging Curtains Conjecture, stated below, points to an important role for the Tur\'an graphs in this problem. 
The Tur\'an graph $T_r(n)$ is the unique graph on $n$ vertices having the maximum number of edges subject to not containing a complete subgraph on $r+1$ vertices. According to Tur\'an's
Theorem, it is the complete $r$-partite graph with parts of nearly equal sizes $\lfloor n/r\rfloor$ or $\lceil n/r\rceil$ (that is, the vertices are partitioned into parts of these sizes, and the graph has all edges joining vertices in different parts\cite{turan}).

In 1997, Masanobu Kaneko defined \emph{poly-Bernoulli numbers}, which bear much the same relation to polylogarithms as Berunoulli numbers do to logarithms. In 2008, Chet Brewbaker described a counting problem whose solution can be identified with the poly-Bernoulli numbers with negative index, the \emph{lonesum matrices}.

Our research was inspired by work on the maximum number of acyclic orientations of graphs with prescribed numbers $n$ and $m$ of vertices and edges respectively, in particular on the ``Hanging Curtains Conjecture'',  which asserts that Tur\'an graphs provide maximum points to the function $m\mapsto a(n,m)$ is concave for $n$ fixed and $m$ between the numbers of edges of successive points. This paper provides some interesting parallel results for graphs restircted to being $r$-partite, and limited evidence for the conjecture.

This paper explores the counting space of acyclic orientations of complete $r$-partite graphs, focusing mainly on  bipartite graphs. There are many strands to the enquiry some more developed than others.
After introducing notation and known results in the preliminary section we go on to develop theoretical properties of the numbers of acyclic orientations of complete bipartite graphs, in section 3. 
This turns out to also be given by \emph{poly-Bernoulli numbers} with negative index~\cite{kaneko}. 
Interestingly the connection was discovered using the On-Line Encyclopedia of Integer Sequences \cite{oeis}.
The connection is proved in two ways in this paper, directly from the formula for the number of acyclic orientation of complete bipartire graphs developed in Section 3.1, and indirectly through the bijection with lonesome marices. In Section 3.4 we show that the sequence  $a(K_{k,n-k})$ for $1\leq k\leq n$ is unimodal, proving that among complete bipartite graphs, the Turán graph $T_2(n)$ maximizes the number of acyclic orientations. 

In Section 4 we explore the space of the number of acyclic orientations of $r$-partite graphs beyond the bipartite case.
The initial exploration, in Section 4.1, reveals the connections between complete bipartite graphs and the number of acyclic orientations of graphs in the neithbourhood of these, in particular revealing their local optimality. 
For $r=3$ we derive a formula for the numbers of acyclic orientations of complete tripartite graphs as well, in Section 4.2.
We also determine the numbers of acyclic orientations of various graphs in the neighbourhood of complete bipartite graphs, and verify that none of these refute the conjecture that the complete bipartite graphs with nearly equal parts are optimal. 

For illustrative purposes, numerical values of complete bipartite graphs and their neighbours, with a removed or added edge or both, is presented in Section 5. This provides some evidence for the Hanging Curtain Conjecture and leads to further speculative conjectures. We provide a range of computational and theoretical results designed to stimulate further research.

\section{Preliminary results}
In this section we introduce notation and definitions used throughout the paper, and provide known results which we need to state and prove our results.\par 

Let $G$ and $G'$ be graphs. Define $G+G'$ to be the graph $G$ and $G'$ and all the possible edges between them. We define $K_{n_1,n_2,\dots,n_r}$ to be the \emph{complete $r$-partite
graph on $n_1 + n_2+\dots +n_r$ vertices} to be the graph whose vertices are partitioned
into sets of sizes $n_i$ for $1\leq i\leq r$, having all possible edges between these
$r$ sets and none within them. Denote the vertices in the part of size $n_i$ by $V_{i}$. Let $N_n$ denote the graph with $n$ vertices and no edges.
Then $$K_{n_1,n_2,\dots,n_r}=K_{n_1}+K_{n_2,\dots,n_r}=\dots=K_{n_1}+K_{n_2}+\dots+K_{n_r.}$$\par  

We first mention a result from Stanley \cite[p.172]{stanley}, which we will use repeatedly throughout. 
\begin{lemma}[Deletion-Contraction Formula]\label{delcont} 
Let $G$ be a graph, and $e$ an edge in $G.$ Then
$$a(G)=a(G-e)+a(G/e).$$
\end{lemma}

We will also use the following properties of the \emph{Stirling number
of the second kind}, denoted by $S(n,k),$ which counts the number of ways to partition a set of $n$
objects into $k$ non-empty subsets. 
Note that $S(n,k)=0$ for $k\leq 0$ and $k\geq n+1.$ 
This satisfies the following iterative relationship derived from the observation that $n^{th}$ object can either be in its own subset or join one of $l$ subsets of the other $n-1$ objects.
\begin{lemma}[Iterative Relationship]
\begin{equation}\label{homokora}
    S(n,l)=S(n-1,l-1)+lS(n-1,l).
\end{equation}
\end{lemma}
Moreover, the following log-concavity property holds 
\begin{lemma}[Log-concave property]
\begin{equation}\label{logconcav}
    S(n-1,k-1)^2>S(n-1,k-2)S(n-1,k).
\end{equation}
\end{lemma}

\section{The number of acyclic orientations of compete bipartite graphs}
The formula for the number of acyclic orientations of a bipartite graph $K_{n_1,n_2}$ is presented in this section. While it might initially appear to be obscure, we show that it is actually the poly-Bernoulli number in the variables $n_1$ and $n_2$, as defined by Kaneko, for which another combinatorial interpretation was found
by Brewbaker. We also demonstrate the equivalence between the number of acyclic orientations of a complete bipartite graph and the number of its Lonesum matrices, thus providing an alternative proof. We then use the formula which we established to describe further characteristics of poly-Bernoulli numbers, along the lines of the Hanging Curtain Conjecture for the maximum number of acyclic orientations of all graphs.

\subsection{The Formula}
In this section we derive a formula for the number of acyclic orientations of a complete bipartite
graph, $K_{n_1,n_2}$ ready for comparison with the poly-Bernoulli numbers in the next section. 
\begin{theorem}
\label{thm:aosknn}
The number of acyclic orientations of the complete bipartite graph
$K_{n_1,n_2}$ is
\[a(K_{n_1,n_2})=\sum_{j=1}^{\min\{n_1+1,n_2+1\}}(j-1)!^2S(n_1+1,j)S(n_2+1,j),\] 
where $S$ denotes Stirling numbers of the second kind.
\end{theorem}

\begin{proof}
Let $A$ and $B$ denote the two bipartite vertex sets. Now any acyclic orientation of the graph can be obtained by ordering the vertices and making the edges point from smaller to greater. If we do this, we will have alternating $A$ and $B$ intervals within the vertex ordering; the ordering within each interval is irrelevant in identifying the orientation since there are no edges involved, but the ordering of the intervals themselves matters.
In terms of structure for a given orientation, call two points $a_1,a_2\in A$ \emph{equivalent} if the orientations of $\{a_1,b\}$ and $\{a_2,b\}$ are the same for all $b\in B$; in other words $a_1$ and $a_2$ are not separated by a vertex in $B$ in any vertex ordering which gives rise to the orientation. Thus, an orientation is uniquely defined by inter-leaving intervals of equivalent vertices in A, respectively B. 

It is left to count alternating intervals. To get around the problem that the first interval in the ordering might be in either $A$ or $B$, and similarly for the last interval, we use the following
device. When the first equivalence vertex subset of the orientation is in $B$ add a dummy set ${a_0}$ of $A $ before it. Similarly if the last equivalence subset is in $A$ insert a dummy subset ${b_0}$ of $B$ after it. 
The mechanism for achieving this construction across all orientations is to take equal numbers of subsets of $A \cup {a_0}$ and $B \cup {b_0}$, $j$ say, and order the subsets so that
\begin{itemize}\itemsep0pt
\item the interval containing $a_0$ is first;
\item followed by alternating $A, B$ bipartite intervals (if any);
\item the interval containing $b_0$ is last.
\end{itemize}
Finally, delete the dummy vertices $a_0$ and $b_0$.

For a given value of $j,$ partitioning of the vertex sets can be done in $S(n_1+1,j)S(n_2+1,j)$ ways, and the inter-leaving of intervals in $(j-1)!^2$ ways. Summing over $j$ gives the total number of ordered intervals and hence of acyclic orientations as claimed.
\end{proof}

\subsection{Poly-Bernoulli numbers}

This is only a very brief introduction to the poly-Bernoulli numbers, which were
introduced by Masanobu Kaneko~\cite{kaneko} in 1997. Kaneko gave the following
definitions. Let
\[\mathrm{Li}_k(z)=\sum_{i=1}^\infty\frac{z^i}{i^k},\]
and let
\[\frac{\mathrm{Li}_k(1-\mathrm{e}^{-x})}{1-\mathrm{e}^{-x}}=
\sum_{n=0}^\infty B_n^{(k)}\frac{x^n}{n!}.\]
The numbers $B_n^{(k)}$ are the \emph{poly-Bernoulli numbers of order~$k$}. Kaneko gave a couple of nice formulae for the poly-Bernoulli numbers of
negative order, of which one is relevant here.

\begin{theorem}[Kaneko]
\[B_n^{(-m)}=\sum_{j=0}^{\min(n,m)}(j!)^2S(n+1,j+1)S(m+1,j+1).\]
\end{theorem}

This formula has the (entirely non-obvious) corollary that these numbers have a symmetry property: $B_n^{(-m)}=B_m^{(-n)}$ for all non-negative integers $n$ and $m$. Kaneko's Theorem together with Theorem \ref{thm:aosknn} gives
the following result.

\begin{theorem}\label{thm3.1}
The number of acyclic orientations of $K_{n_1,n_2}$ is the poly-Bernoulli number $B_{n_1}^{(-n_2)}=B_{n_2}^{(-n_1)}$.
\end{theorem}

\subsection{Lonesum matrices}

Another combinatorial interpretation, in terms of lonesum matrices, was given by Chad Brewbaker~\cite{cb}
in 2008. A zero-one matrix is a \emph{lonesum matrix} if it is uniquely
determined by its row and column sums. Clearly a lonesum matrix cannot contain either
$\begin{pmatrix}
  1 & 0\\ 
  0 & 1
\end{pmatrix}$ or
$\begin{pmatrix}
  0 & 1\\ 
  1 & 0
\end{pmatrix}$  as a submatrix (in not necessarily
consecutive rows or columns, since if one such submatrix occurred it could be
flipped into the other without changing the row and column sums.)
Ryser~\cite{ryser} showed that, conversely, a matrix containing neither of
these is a lonesum matrix. 
\begin{theorem}[Brewbaker]
The number of $n_1\times n_2$ lonesum matrices is given
by the poly-Bernoulli number $B_{n_1}^{(-n_2)}.$
\end{theorem}
We give a simple argument why this number is equal
to the number of acyclic orientations of $K_{n_1,n_2}$ and hence, given Theorem \ref{thm3.1}, provide an alternative proof to Brewbaker's result.

\begin{theorem}
The number of acyclic orientations of the bipartite graph $K_{n_1,n_2}$  is equal to the number of $n_1\times n_2$ lonesum matrices. 
\end{theorem}
\begin{proof}

In one direction, number the vertices
in the bipartite blocks from $1$ to $n_1$ (in $A$) and from $1$ to $n_2$
(in $B$). Now given an orientation of the graph, we can describe it by a matrix
whose $(i,j)$ entry is $1$ if the edge from vertex $i$ of $A$ to vertex $j$
of $B$ goes in the direction from $A$ to $B$, and $0$ otherwise. The two
forbidden submatrices for lonesum matrices correspond to directed $4$-cycles;
so any acyclic orientation gives us a lonesum matrix.

Conversely, if an orientation of a complete bipartite graph contains no
directed $4$-cycles, then it contains no directed cycles at all. For suppose that there are no directed $4$-cycles, but there is a
directed $2k$-cycle for $k\geq 3,$ $(a_1,b_1,a_2,b_2,\ldots,a_k,b_k,a_1)$. Then the edge between
$a_1$ and $b_{k-1}$ must be directed from $b_{k-1}$ to $a_1$, since otherwise there
would be a $4$-cycle $(a_1,b_{k-1},a_k,b_k,a_1)$. But then we have a shorter
directed cycle, $(a_1,b_1,a_2,\ldots,b_{k-1},a_1)$. Continuing this shortening
process, we would eventually arrive at a directed $4$-cycle, a contradiction.
(This simply says that the cycle space of the complete bipartite graph is
generated by $4$-cycles.)
\end{proof}

\subsection{Unimodality of $B_{k}^{(-n+k)}$ and of $a(K_{k,n-k})$}
In this section we will show unimodality of the sequence $a(K_{k,n-k})$ for $k=1,2,\dots,n-1.$ 
A sequence $x_0,\dots,x_{n}$ is \emph{unimodal} if there is some $0\leq k\leq n$ such that $x_0\leq x_1\leq \dots,\leq x_k\geq x_{k+1}\geq \dots \geq x_{n}.$ 

\begin{theorem}\label{unimodal}
 The number of acyclic orientations over the class of all complete bipartite graphs with the same number, $n,$ of vertices, is maximised by a graph $K_{n_1,n_2}$, for which $|n_1-n_2|\leq 1,$ namely the Turán graph $T_2(n)=K(\lfloor n/2\rfloor,\lceil n/2 \rceil).$
 \\Moreover, the sequence $a(K_{k,n-k})$ for $k=1,2,..,n-1$ is unimodal.
\end{theorem}
By Theorem \ref{thm3.1} we then get the following corollary.
\begin{cor}
The poly-Bernoulli numbers are unimodal with maximum value $B_{n_1}^{(-{n_2})}$ for $|n_1-n_2|\leq 1$.
\end{cor}

To prove Theorem \ref{unimodal}, we first establish progressively sophisticated properties of the Stirling numbers of the Second Kind building upon known results mentioned in the Preliminary Section.

\begin{lemma}\label{lemma:10}
$$\frac{S(n-1,k-1)}{S(n-1,k)}>\frac{S(n,k-1)}{S(n,k)}-\frac{S(n-1,k-1)}{S(n,k)}$$
$$>\frac{S(n,k-1)}{S(n,k)}$$
\end{lemma}
\begin{proof} Applying the primary property of $S(n,l),$ Equation \ref{homokora}, for $n$ with $l=k$ and $n-1$ with $l=k-1,$ respectively, gives
\begin{equation}\label{A}
    S(n-1,k-1)S(n,k)=S(n-1,k-1)^2+kS(n-1,k)S(n-1,k-1)
\end{equation}
and
\begin{equation}\label{B}
    S(n-1,k)S(n,k-1)=S(n-1,k-2)S(n-1,k)+(k-1)S(n-1,k)S(n-1,k-1).
\end{equation}
Subtracting expression \ref{B} from \ref{A} and applying the log-concavity property of the Stirling numbers of the second kind, Equation \ref{logconcav}, gives 
$$S(n-1,k-1)S(n,k)-S(n-1,k)S(n,k-1)>S(n-1,k)S(n-1,k-1)$$ 
and the result follows. 
\end{proof}

\begin{lemma}\label{thm:55}
Let $n_1$ and $n_2$ be natural numbers such that $n_1\leq n_2$. Then, for $m\leq n_1 -1,$
$$S(n_1,k)S(n_2,k)>S(n_1-m,k)S(n_2+m,k).$$ 
\end{lemma}

\begin{proof}
We first prove the result for $m=1$, namely, that
$$S(n_1,k)S(n_2,k)>S(n_1-1,k)S(n_2+1,k),$$
or, after rearranging,  
\begin{equation}\label{alpha}
    \frac{S(n_1,k)}{S(n_1-1,k)}>\frac{S(n_2+1,k)}{S(n_2,k)}.
\end{equation}
Applying property \ref{homokora} to the numerators gives the following expressions for the two terms $$\frac{S(n_1,k)}{S(n_1-1,k)}=\frac{S(n_1-1,k-1)}{S(n_1-1,k)}+k$$ and 
$$\frac{S(n_2+1,k)}{S(n_2,k)}=\frac{S(n_2,k-1)}{S(n_2,k)}+k.$$ Then, since $n_2\geq n_1-1,$ repeated application of Lemma \ref{lemma:10} applied to the term in the right hand side of the above expressions results in inequality \ref{alpha} as required.
The general result follows by applying the inequality $m$ times for decreasing values of $n_1$ (and increasing values of $n_2$.)
\end{proof}

We are now ready to prove Theorem \ref{unimodal}.

\begin{proof}[Proof of Theorem \ref{unimodal}] We begin by showing that $K_{n_1,n_2}$, for $|n_1-n_2|\leq 1$, maximises the number of acyclic orientations over the class of all complete bipartite graphs with the same number of vertices. 
Let $K_{a,b}$ be an arbitrary bipartite graph on $n=n_1+n_2$ vertices differing from $K_{n_1,n_2}$, i.e. $|b-a|\geq 2$. We will now show that $a(K_{a,b})<a(K_{n_1,n_2})$. Let $a<b$ and $n_1<n_2$. Then clearly also $a<n_1$, so we can apply Lemma \ref{thm:55} alongside Theorem \ref{thm:aosknn} to get $$a(K_{a,b}) \\ = \sum_{k=1}^{a+1} ((k-1)!)^2S(a+1,k)S(b+1,k) < \sum_{k=1}^{a+1} ((k-1)!)^2S(n_1+1,k)S(n_2+1,k).$$ 
Also, since $a+1 < n_1+1,$ applying Theorem \ref{thm:aosknn} to $a(K_{n_1,n_2}),$ 
$$
a(K_{n_1,n_1})=
\sum_{k=1}^{n_1+1}((k-1)!)^2S(n_1+1,k)S(n_2+1,k) >
\sum_{k=1}^{a+1} ((k-1)!)^2S(n_1+1,k)S(n_2+1,k).
$$
Thus, $$a(K_{a,b})<a(K_{n_1,n_2})$$
and the result follows. \par 
Note that complete bipartite graphs are symmetric, $K_{k,n-k}$ and $K_{n-k,n-(n-k)}$ are in fact isomorphic. Hence the sequence $a(K_{k,n-k})$ for $k=1,2,..,n$ is symmetric. Now unimodality follows from above.
\end{proof}

We conjecture that the unimodularity property also holds for the class of complete $r$-partite graphs, for any value $r$, and in particular that the Tur\'an graphs play a special role.

\begin{conjecture}\label{conj:rpart}
For the class of complete $r$-partite graphs on $n$ vertices, the number of acyclic orientations has a concavity property with respect to the number of edges, and is maximised by the Tur\'an graph, $T_{(r,n)}.$ 
\end{conjecture}

We might speculate more generally.

 \begin{conjecture}
Let $G_1,$ $G_2$ and $H$ be graphs. If $a(G_1)\leq{a(G_2)}$ then $a(G_1+H)\leq a(G_2+H)$. 
\end{conjecture}

The concavity property of conjecture \ref{conj:rpart}, would follow directly from this more general conjecture. To see this consider the following iterative step between complete $r$-partite graphs on $n$ vertices. 
While any pair of vertex set partitions have sizes differing from $\lceil n/r \rceil$ and $\lfloor n/r \rfloor,$ apply Theorem \ref{unimodal} to the corresponding bipartite subgraph, $G_1,$ to find a substitute bipartite graph. $G_2$ with more equal partition size without increasing the number of acyclic orientations.

\section{Exploration of acyclic orientation counts on $r$-partite graphs}
We explore the solution space of the number of acyclic orientations close to complete bipartite graphs, differing by 1 or 2 edges.
In these proofs we identify a graph of particular interest related to both $K_{n_1,n_2}$ and $K_{n_1,n_2-1}$
(or $K_{n_1-1,n_2}$) which also provides the simplest step towards a complete tripartite graph, $K_{n_1-1,n_2-1,1}.$ 
In fact, this theorem generalises to $r$-partite graphs.
To round off the section we derive a formula for the number of acyclic orientation of complete tri-partite graphs and offer  an algorithmic approach for $r$-partite ones. \par 

\subsection{Local optimality of complete bipartite graphs}
We will now establish the local maximality with respect to the number of its acyclic orientations  of a complete bipartite graph compared to its closest neighbours with the same number of vertices and edges. Note that the results and proofs in this section naturally generalise to further dimensions of the $r$-partite graphs.

\begin{theorem}\label{thm:3}
Let $e$ and $e_1$ as in Theorem \ref{thm:aosknnpluse}. Then for $n_1\geq 2,$ $n_2\geq 2$
\begin{align*}
    a(K_{n_{1},n_{2}})&> a(K_{n_1,n_2}+e_{1}-e).
\end{align*} 
\end{theorem}

In preparation for the proof of \ref{thm:3} we start by establishing connections between compete bipartite and tripartite graphs in the following lemma, and then using it to establish the relationship between the number of acyclic orientations of the various relevant graphs. These constructions also provide the basis for the formula of the number of acyclic orientations of tripartite graphs in the next section.
\begin{lemma}\label{lemmatrip}
Let $K_{n_1,n_2} = (V_1 \cup V_2,E)$ and $w$ denote the vertex in the third partition of $K_{n_1,n_2,1}.$ Then
\begin{align}
  K_{n_1,n_2,1}&=K_{n_1+1,n_2+1}/e\nonumber\\
                &=K_{n_1,n_2}+N_1\nonumber\\
                &=K_{n_1,n_2+1}\cup E_2 \nonumber
\end{align}
where $E_2=\{(v,w))\vert v\in V_2\}.$ 
\end{lemma}

\begin{proof}
Concatenating an edge $e$ in $K_{n_1+1,n_2+1}$ produces a vertex connected to all other vertices in the graph, while leaving the remaining subgraph on $n_1+n_2$ vertices, $K_{n_1,n_2},$ untouched. This graph may also be described as adding the empty graph on 1 vertex to $K_{n_1,n_2},$ i.e. $K_{n_1,n_2}+N_1,$ which is by definition in fact $K_{n_1,n_2,1}.$ \par
Observe that $K_{n_1,n_2,1}$ is of the form $G(V_1 \cup V_2 \cup w,E\cup E_1\cup E_2)$ where $E_i=\{(v_i,w))\vert v_i\in V_i\}$ for $i=1,2,$ and that in $K_{n_1,n_2}+N_1$ the single vertex, $w,$ of $N_1$ in the third partition is connected to all vertices in the first partition.
Thus, $K_{n_1,n_2+1}$ can be viewed as a subgraph of $K_{n_1,n_2,1},$ as the edges of the set $E_2$ are the only things missing.
\end{proof}

We now explore the graph space around the complete bipartite graph. We first add an edge, then we move an edge while minimising disruption of vertex degrees. We then compare the number of acyclic orientations of these closely related graphs to that of the complete bipartite graph.

\begin{lemma}
\label{thm:aosknnpluse}
The number of acyclic orientations of a bipartite graph, when $n_1\geq 2,$ is related to that of neighbouring graphs differing by at most 2 edges as follows: 
\begin{align}
    a(K_{n_1,n_2}+e_1)&=a(K_{n_1,n_2})+a(K_{n_1-1,n_2})\label{x}\\
    a(K_{n_1,n_2}+e_1-e)&=a(K_{n_1,n_2}-e)+a(K_{n_1-1,n_2})\label{y}\\
     a(K_{n_1,n_2})&=a(K_{n_1,n_2}-e)+a(K_{n_1-1,n_2-1,1})\label{z},
\end{align}
where $e$ is an edge in $K_{n_1,n_2}$ and $e_1$ an edge joining its endpoint in $V_{n_1}$ to another vertex in $V_{n_1}$.
\end{lemma}

\begin{proof}
We apply the Deletion-Contraction formula in Lemma \ref{delcont} to each of the terms on the left hand side of each of the three expressions. It is left to identify the graphs in the second term on the right hand side with the corresponding concatenated graph for each of the three expressions. \par 
Observe that $(K_{n_1,n_2}+e_1)/e_1$ concatenates 2 vertices in $V_{n_1}$ into 1, resulting in $K_{n_1-1,n_2},$ as required.\par 
Similarly, $(K_{n_1,n_2}+e_1-e)/e_1$ concatenates 2 vertices in $V_{n_1}$, one of which is adjacent to all vertices of $V_{n_2}$, and so the concatenated vertex is thus also, and the resulting graph is $K_{n_1-1,n_2}.$ \par 
For the third expression, the last term $K_{n_1,n_2}/e$ is shown in Lemma \ref{lemmatrip} to be the tripartite graph $K_{n_1,n_2,1},$ completing the proof.
\end{proof}

\begin{cor}

\label{cor:3}
Let $e$ and $e_1$ as in Theorem \ref{thm:aosknnpluse}. Then for $n_1\geq 2,$ $n_2\geq 2$ 
    $a(K_{n_{1},n_{2}})=a(K_{n_1,n_2}+e_1-e)+a(K_{n_1-1,n_2-1,1})-a(K_{n_1-1,n_2})$.
\end{cor}

 \begin{proof}
 Combining equation \ref{y} and \ref{z} gives the following relationship
 $a(K_{n_{1},n_{2}})=a(K_{n_1,n_2}+e_1-e)+a(K_{n_1-1,n_2-1,1}-a(K_{n_1-1,n_2}).$
 \end{proof}
 
 We are now ready to prove Theorem \ref{thm:aosknnpluse}.
 \begin{proof}
From corollary \ref{cor:aobipar=}, it is left to show that $a(K_{n_1-1,n_2-1,1}) > a(K_{n_1-1,n_2})$ to establish the inequality. 
Applying Lemma \ref{lemmatrip} followed by the Deletion-Contraction Lemma \ref{delcont} applied iteratively to remove all but one of the edges of $E_2,$ gives \begin{align*}
    a(K_{n_1-1,n_2-1,1})=a(K_{n_1-1,n_2}+E_2)&\geq a(K_{n_1-1,n_2}+e_2)
\end{align*} 
where $e_2$ is an arbitrary edge between vertices in $V_{n_2}.$ Now applying equation \ref{x} with the roles of $V_1$ and $V_2$ reversed gives $a(K_{n_1-1,n_2}+e_2)=a(K_{n_1-1,n_2})+a(K_{n_1-1,n_2-1}).$ The required inequality follows as $a(K_{n_1-1,n_2-1})\geq 1,$ since $n_1-1\geq1$ and $n_2-1\geq1.$ 
 \end{proof}
 
\begin{cor}\label{cor:aobipar=}
Let $e$ and $e_1$ as in Theorem \ref{thm:aosknnpluse}. Then for $n_1\geq 2,$ $n_2\geq 2$
\begin{align*}
    a(K_{n_{1},n_{2}})&=a(K_{n_1,n_2}+e_1-e)+a(K_{n_1-1,n_2-1,1})-a(K_{n_1-1,n_2})\\
\end{align*} 
\end{cor}


\subsection{A formula for a specific complete tripartite graph}
The following theorem gives a formula for the number of acyclic orientations of complete tripartite graphs of the form $K_{n_1-1,n_2-1,1}$ and, in view of Theorems \ref{thm:aosknnpluse}, equation \ref{z} and Theorem \ref{thm:aosknn}, also for $a(K_{n_1,n_2}-e).$

\begin{theorem}
\label{thm:aosknnminuse} The number of acyclic orientations of the complete bipartite graph $K_{n_1-1,n_2-1,1}$ is
$$a(K_{n_1-1,n_2-1,1})=\frac{X}{2},$$ where
\begin{eqnarray*}
X=1&&+\sum_{i=2}^{\min\{n_1,n_2\}+1} \quad((i-2)!)^2\\
&&[ (2i-3)S(n_1+1,i)S(n_2+1,i)
- (i-2)(S(n_1+1,i)S(n_2,i)\\
&&+S(n_1,i)S(n_2+1,i))
-       S(n_1,i)  S(n_2,i)\,].
\end{eqnarray*}
\end{theorem}

\begin{proof}

Now, consider the acyclic orientations of $K_{n_1,n_2}-e$. They extend to either 2 or 1 acyclic orientations of $K_{n_1,n_2}$ depending upon whether the edge $e$ can be flipped without creating a cycle or not. Let $Y$ and $Z$ be the number of acyclic orientations of each type, respectively, in $K_{n_1,n_2}-e$.  Then $a(K_{n_1,n_2}-e) = Y + Z$, $a(K_{n_1,n_2}) = 2Y + Z$, and $X = 2Y$, so
 $a(K_{n_1,n_2}) - a(K_{n_1,n_2}-e) = Y$.
Now from Lemma \ref{lemmatrip} and the Deletion-Contraction formula,
$a(K_{n_1-1,n_2-1 ,1}) = a(K_{n_1,n_2}/e) = a(K_{n_1,n_2}) - a(K_{n_1,n_2}-e)= Y$
It thus remains to verify the formula for $X=2Y$ given in the statement of the theorem.

Then $X=2Y$. Thus, it remains to count the number of acyclic orientations of $K_{n_1,n_2}$ in which the edge $(a,b)$ can be flipped.
Let $X$ be the number of acyclic orientations of $K_{n_1,n_2}$ which remain acyclic
when a given edge  $e=\{a,b\}$ is flipped. 
(This number clearly does not depend on the chosen edge.)

We use the description of the number of acyclic orientations of $K_{n_1,n_2}$ as ordered intervals of $A$ and $B.$
Observe that the edge $e=(a,b)$ can be flipped in an orientation without creating a cycle if and only if
the part of the partition of $B$ containing $b$ immediately precedes or
follows the part of the partition of $A$ containing $a,$ in the corresponding vertex ordering. 
(If a part of $B,$ containing a vertex $b'$ say, and a part of $A,$ containing a vertex $a',$ intervene,
then we have arcs $(a,b')$, $(b',a')$ and $(a',b)$, so the arc
$(a,b)$ is forced. Similarly in the other case.)

To establish the count $X$ we follow the construction in the proof of Theorem~\ref{thm:aosknn}.

The case $j=1$ applies to the single orientation in which all edges are directed from $A$ to $B$, and $(a,b)$ can be flipped. 
So this contributes $1$ to the sum making up $X.$
Suppose that $j\geq 2$. We distinguish four cases, according as $a_0$ and $a$ are or are not in the same part, and similarly for $b_0$ and $b$. 
Of the $S(n_1+1,j)$ partitions of $A\cup\{a_0\}$, $S(n_1,j)$ have $a_0$ and $a$ in the same part: this is found by regarding $a_0$ and $a$ as the same element, partitioning the resulting set of size $n_1$, and then separating them again.

\subparagraph{Case 1} $a_0$ and $a$ in the same part, $b_0$ and  $b,$ in
the same part. 
Since $j>1$, the parts containing $a_0$ and $b_0,$ and hence the parts containing $a$ and $b,$ are not
consecutive, so the contribution from this case is $0$.

\subparagraph{Case 2} $a_0$ and $a$ in the same part, $b_0$ and $b$ not.
There are $S(n_1,j)(S(n_2+1,j)-S(n_2,j))$ pairs of partitions with this
property. Now the part containing $b$ must come immediately after the
part containing $a$, so there are only $(j-2)!$ orderings of the parts
of $B$, while still $(j-1)!$ for the parts of $A$.

\subparagraph{Case 3} $b_0$ and $b$ in the same part, $a_0$ and $a$ not.
This case is the same as Case 2, with $n_1$ and $n_2$ interchanged.

\subparagraph{Case 4} $a_0$ and $a$ in different parts, $b_0$ and $b$ in different parts. 
There are $(S(n_1+1,j)-S(n_1,j))(S(n_2+1,j)-S(n_2,j))$
such pairs of partitions. 
Now the parts containing $a$ and $b$ must be adjacent, so must occur as $(3,2),(3,4),(5,4),\ldots,$ or $(2j-1,2j-2)$ in the ordering of parts: there are $(2j-3)$ possibilities. Once one possibility has been chosen, the position of two intervals containing $a$ and $b$ are fixed, leaving $((j-2)!)^2$ possible orderings of the other intervals.

\medskip
Combining all of the above terms and rearranging, gives the value of $X,$ completing the proof.
\end{proof}

\subsection{An algorithm to compute the number of acyclic orientations of $r$-partite graphs }

By \cite{stanley}, we can compute the number of acyclic orientations of a graph by substituting $-1$ into its chromatic polynomial. 

In this section we produce an iterative algorithm to find the chromatic polynomial of many families of graphs, including all complete $r$-partite graphs.\par  
Let $P_G(q)$ denote the chromatic polynomial of $G$. Moreover, let $P_G^*(q)$
(for a positive integer $t$) be the number of colourings of $G$ with $q$
colours, all of which are used. Then the Inclusion--Exclusion Principle gives
\[P_G^*(q)=\sum_{s=1}^q(-1)^{q-s}{q\choose s}P_G(s).\]
We can include $s=0$ in the sum since $P_G(0)=0$.

Note that $P_G^*(q)$ is not a polynomial in $q$, since it is zero for all $q$
greater than the number of vertices of $G$. Moreover, it does not make sense
to substitute a negative number for $q$ in this formula.

Let $G+H$ denote the sum graph, the disjoint union of $G$ and $H$ with all
edges between $G$ and $H$. Then
\[P_{G+H}(q)=\sum_{r=0}^q{q\choose r}P_G^*(r)P_H(q-r).\]
This is because the two summands must use disjoint sets of colours; so if
a prescribed set of $r$ colours is used for $G$, all of them used, then $H$
must be coloured with the remaining $q-r$ colours.

This formula has two disadvantages. First, it is not symmetric between $G$
and $H$. Second, although $P_{G+H}(q)$ is a polynomial in $q$, we cannot
just substitute $q=-1$ in the formula to obtain the number of acyclic
orientations, since we have a sum from $0$ to $q$. What follows is an
attempt to remedy these defects.

\begin{theorem}
\[P_{G+H}(q)=\sum_{s,t}^{s+t \leq q}(-1)^{q-s-t}\frac{q!}{s!t!(q-s-t)!}P_G(s)P_H(t),\]
where the summation is over pairs $(s.t)$ of non-negative integers with sum at
most~$q$.
\end{theorem}

\begin{proof}

\begin{eqnarray*}
P_{G+H}(q) &=& \sum_{r=0}^q{q\choose r}P_G^*(r)P_H(q-r)\\
&=& \sum_{r=0}^q\sum_{s=0}^r(-1)^{r-s}{q\choose r}{r\choose s}P_G(s)P_H(q-r).
\end{eqnarray*}

The sum is over all pairs $(r,s)$ with $0\le s\le r\le q$.
Now put $t=q-r$; some manipulation now gives the result, using
the fact that
\[{q\choose r}{r\choose s}=\frac{q!}{s!(r-s)!(q-r)!}.\]
\end{proof}
This expression is now symmetric between $G$ and $H$. However, we still cannot
substitute $q=-1$. Rather, we should regard the expression as an algorithm
for calculating the chromatic polynomial of $G+H$ explicitly; then we can
substitute $q=-1$ into this polynomial.

\begin{cor}
$P_{G+H}(q)=(-1)^q \sum {q\choose s+1}\frac{(-1)^s}{s!}P_G(s)\frac{(-1)^t}{t!}P_H(t).$
\end{cor}
\begin{proof}

Writing the trinomial coefficient as
\[{q\choose s+t}{s+t\choose t},\]
the dependence on $q$ is entirely in the factors $(-1)^{q-s-t}{q\choose s+t}$.
\end{proof}
The complete bipartite graph is the sum of two edgeless graphs; so we get a
formula for its chromatic polynomial:
\begin{cor}
The chromatic polynomial for the bipartite graph as a sum
of two edgeless graphs is
\[P_{K_{m,n}}(q)=
(-1)^q\sum_{s,t}^{s+t\leq q}{(q)_{s+t}}\frac{(-1)^ss^m}{s!}\frac{(-1)^tt^n}{t!},\]
where $(q)_{s+t}$ denotes the falling factorial $q(q-1)\cdots(q-s-t+1)$.
\end{cor}
\section{Computational Evidence}

It is instructive to consider numerical results in parallel to provable properties of $a(K_{n_1,n_2})$ and
neighbouring values for both testing and inspiring conjectures maximality and the shape of the
solution space.
We performed two sets of computations, one for all small values of vertex sets $n_1$ and $n_2$ up to
7, and the other for values of $n_1$ and $n_2$ both equal to 10, 100 or 1000 (or 2000). The
computational intensity of the work limited fuller evaluations. In both cases we explore the effects
of adding and/or deleting an edge to understand the nature of the solution space, with respect to
the number of cyclic orientations, around the complete bipartite graphs. The results illustrate the
concavity property of the Hanging Curtain Conjecture.
The computations are from the fourth author's PhD thesis\cite{robert}.

It is instructive to view the numerical values of the number of acyclic orientations of bipartite graphs $K_{n_1,n_2}$. For $n_1$ between $2$ and $7$ Table \ref{tab:aosofknn} gives the number of acyclic orientations of the complete bipartite
graphs.  \par When $n_1=1$, the graph $K_{1,n}$ is a tree, and we have $a(K_{1,n})=2^n$. Similarly we have $a(K_{2,n}+e_1)=2\cdot3^n$. This is because the graph $K_{2,n}+e_1$
consists of $n$ triangles sharing a common edge $e_1$, there are two ways to
orient the edge $e_1$, and then three ways to choose the orientations of the remaining edges of each triangle
to avoid a cycle. Putting these two results together in Theorem \ref{thm:aosknnpluse} gives us $a(K_{2,n})=2\cdot3^n-2^n$.  
The rest of the values in the tables have been checked by calculating the chromatic polynomial
of the graph.  (Recall that a theorem of Stanley~\cite{stanley} asserts that the number of
acyclic orientations of an $n$-vertex graph $G$ is $(-1)^nP_G(-1)$, where
$P_G$ is the chromatic polynomial of $G$.)

\begin{table}[htbp]
\[\begin{array}{|c||r|r|r|r|r|r|}
\hline
n_1\setminus n_2 & \multicolumn{1}{c|}{2} & \multicolumn{1}{c|}{3} &
\multicolumn{1}{c|}{4} &  \multicolumn{1}{c|}{5} &  \multicolumn{1}{c|}{6} &
 \multicolumn{1}{c|}{7} \\
\hline
\hline
2 & 14 & 46 & 146 & 454 & 1394 & 4246 \\
\hline
3 & & 230 & 1066 & 4718 & 20266 & 85310 \\
\hline
4 & & & 6902 & 41506 & 237686 & 1315666 \\
\hline
5 & & & & 329462 & 2441314 & 17234438 \\
\hline
6 & & & & & 22934774 & 202229266 \\
\hline
7 & & & & & & 2193664790 \\
\hline
\end{array}\]
\caption{The number of acyclic orientations of $K_{n_1,n_2}$}
\label{tab:aosofknn}
\end{table}

The following table contains the number of acyclic orientations of complete bipartite graphs $K_{n_1,n_2}$ on 8 vertices and $m$ edges. It also contains the maximum number $a_{max}$ of acyclic orientations possible for graphs with $m$ edges. Note that these numbers coincide exactly when $n_1=n_2,$ so the graph is Turán. 

\begin{table}[htbp]
\[\begin{array}{|c|r|r|r|r|}
\hline
n_1& n_2 & m & a_{K_{n_1,n_2}} &a_{max} \\
\hline
\hline
2 & 6 & 12 & 1394 &1920  \\
\hline
3 & 5& 15 & 4718 &5000 \\
\hline
4 & 4& 16& 6902 & 6902 \\
\hline
\end{array}\]
\caption{Maximal acyclic orientations of bipartite graphs for $n=8$}
\label{tab:aosofknnn}
\end{table}



We conjecture that for biparite blocks of equal size, $K_{n,n}+e$ and $K_{n,n}-e$, also maximise
the number of acyclic orientations for graphs with the same number of vertices
and of edges.

Numerical results for values of $n=10$, $100$, $1000$ (and $2000$) provide some insight into the
behaviour of these functions for larger instances.

Related to this conjecture, we observed that the ratio
\[\frac{a(K_{n,n})-a(K_{n,n}-e)}{a(K_{n,n}+e)-a(K_{n,n})}\]
is about 2 within the range of computation; its values for $n=10$, $100$, $1000$
are respectively $1.923534$, $1.992995$, $1.999306$ respectively. The
convergence is quite slow; the computed values appear to be $2-O(n^{-1})$.

We also observed that, within the range of computation, $K_{n,n}-e$ has more
acyclic orientations than $K_{n+1,n-1}$ (these graphs have the same numbers of
vertices and edges). For $n=10$, $100$, $1000$, the ratio
\[\frac{a(K_{n,n})-a(K_{n+1,n-1})}{a(K_{n,n})-a(K_{n,n}-e)}\]
is $1.367903$, $1.596801$, $1.626101$ respectively. It is not so clear how
these ratios behave.

We now consider some related graphs which are ``close to'' complete bipartite
graphs.

Within the range of computation (1,100) $K_{n,n}$ has more acyclic orientations
than $K_{n-1,n+1}+e$, as computed above.  
For  $n=10$, $100$, $1000,$ $2000$ the ratio 
\begin{equation}
\frac{a(K_{n,n})-a(K_{n-1,n+1})}{a(K_{n-1,n+1}+e)-a(K_{n-1,n+1})}
\end{equation}
is $3.921676$, $3.3275021$, $3.2657740$, $3.2623332$. It seems to show some convergence, but not clear. 

Let $A_n(m)$ be the maximum value of $a(G)$ over graphs with $n$ vertices and $m$ edges. 
Now we state the \emph{Hanging Curtains Conjecture}.
\begin{conjecture*}[Hanging Curtains Conjecture]
For a fixed value of $n$ vertices, and value $m$ for which there is a Tur\'an graph, this realises the maximum value $A_n(m)$ for the number of acyclic orientations. This is referred to as a Tur\'an point. The graph of $A_n(m)$ is concave at $m,$ that is $$A_n(m+1)-A_n(m)<A_n(m)-A_n(m-1).$$ Moreover, the graph of $A_n(m)$ for fixed n is roughly convex between successive Tur\'an points.
\end{conjecture*}

In other words, the graph resembles the shape of curtains suspended from the points corresponding to Tur\'an graphs.

\section{Acknowledgements}

Robert Schumacher was funded by EPSRC Grant EP/P504872/1. When joining the project in 2018, Kamilla Rekvényi was funded by a Summer Research Stipend awarded by the School of Mathematics, University of St Andrews.

\bibliographystyle{plain}

\appendix
\section{Data for graphs close to complete bipartite graphs}
Tables \ref{tab:aosknnpluse} and \ref{tab:aosknnminuse} give the number of acyclic orientations of the graphs with an edge added or removed from the complete bipartite graph $K_{n_1,n_2}$, calculated from the formulae in Theorems \ref{thm:aosknn}, \ref{thm:aosknnpluse} and \ref{thm:aosknnminuse}. In Table \ref{tab:aosknnpluse} for
$K_{n_1,n_2}+e_1$, the added edge $e_1$ is in the bipartite block of size $n_1$.
\begin{table}[htbp]
\[\begin{array}{|c|r|r|r|r|r|r|}
\hline
n_1\setminus n_2 & \multicolumn{1}{c|}{2} & \multicolumn{1}{c|}{3} &
\multicolumn{1}{c|}{4} &  \multicolumn{1}{c|}{5} &  \multicolumn{1}{c|}{6} &
 \multicolumn{1}{c|}{7} \\
\hline
2 & 18 & 54 & 162 & 486 & 1458 & 4374 \\
\hline
3 & 60 & 276 & 1212 & 5172 & 21660 & 89556 \\
\hline
4 & 192 & 1296 & 7968 & 46224 & 257952 & 1400976 \\
\hline
5 & 600 & 5784 & 48408 & 370968 & 2679000 & 18550104 \\
\hline
6 & 1848 & 24984 & 279192 & 2770776 & 25376088 & 219463704 \\
\hline
7 & 5640 & 105576 & 1553352 & 19675752 & 225164040 & 2395894056 \\
\hline
\end{array}\]
\caption{The number of acyclic orientations of $K_{n_1,n_2}+e_1$}
\label{tab:aosknnpluse}
\end{table}

\begin{table}[htbp]
\[\begin{array}{|c|r|r|r|r|r|r|}
\hline
n_1\setminus n_2 & \multicolumn{1}{c|}{2} & \multicolumn{1}{c|}{3} &
\multicolumn{1}{c|}{4} &  \multicolumn{1}{c|}{5} &  \multicolumn{1}{c|}{6} &
 \multicolumn{1}{c|}{7} \\
\hline
2 & 8 & 28 & 92 & 292 & 908 & 2788 \\
\hline
3 & & 152 & 736 & 3344 & 14608 & 62192 \\
\hline
4 & & & 5000 & 30952 & 180632 & 1012936 \\
\hline
5 & & & & 253352 & 1915672 & 13715144 \\
\hline
6 & & & & & 18381608 & 164501368 \\
\hline
7 & & & & & & 1812141032 \\
\hline
\end{array}\]
\caption{The number of acyclic orientations of $K_{n_1,n_2}-e$}
\label{tab:aosknnminuse}
\end{table}

\end{document}